\documentclass[11pt]{article}
\usepackage[utf8]{inputenc}

\usepackage{amsthm, amsfonts, amssymb, amsmath,enumitem, natbib, bbm, mathabx}
\usepackage{cases, amsthm}
\usepackage{fullpage}
\usepackage{mathtools}
\mathtoolsset{showonlyrefs}
\usepackage{ wasysym} 
\usepackage{fancyhdr}
\usepackage{graphicx}
\usepackage{bbm}
\usepackage{caption}
\usepackage{subcaption}
\usepackage{algorithm}
\usepackage{algorithmic}
\usepackage{authblk}
 \usepackage{setspace}
\let\Algorithm\algorithm
\renewcommand\algorithm[1][]{\Algorithm[#1]\setstretch{1.6}}

\usepackage[colorlinks=true, pdfstartview=FitV,
linkcolor=blue, citecolor=blue, urlcolor=blue]{hyperref}

\newtheorem{definition}{Definition}[section]
\newtheorem{theorem}[definition]{Theorem}
\newtheorem{proposition}[definition]{Proposition}

\usepackage{xcolor}

\usepackage[top=1in, bottom=1in, left=1in, right=1in]{geometry}
\title{Central Limit Theorems and Approximation Theory: Part II}

\author[1]{Arun Kumar Kuchibhotla}
\affil[1]{Department of Statistics \& Data Science, Carnegie Mellon University}
\date{}

\begin{document}

\maketitle

\begin{abstract}  
In Part I of this article~\citep{banerjee2023central}, we have introduced a new method to bound the difference in expectations of an average of independent random vector and the limiting Gaussian random vector using level sets. In the current article, we further explore this idea using finite sample Edgeworth expansions and also established integral representation theorems. 
\end{abstract}

\section{Introduction}
A sequence $\{W_n\}$ of random variables in a measurable space $\mathcal{W}$ converges weakly to another random variable $W$ if and only if $\mathbb{E}f(W_n)$ converges to $\mathbb{E}f(W)$ for every bounded continuous function $f:\mathcal{W}\to\mathbb{R}$~\citep[Theorem 1.3]{bhattacharya2010normal}. In other words, $$\Delta_f(W_n, W) := |\mathbb{E}[f(W_n)] - \mathbb{E}[f(W)]|,$$ converges to zero for every bounded continuous function $f:\mathcal{W}\to\mathbb{R}$. One can think of this as an asymptotic statement and ask if there is a finite sample version of this result that implies a tractable bound on $\Delta_f(W_n, W)$ that depends on $n$, the distribution of $W_n, W,$ and $f$ and converges to zero as $n\to\infty$. 

The classical central limit theorem implies under a variety of regularity conditions that $W_n = n^{-1/2}\sum_{i=1}^n X_i$ for independent and identically distributed random variables $X_i$ with mean zero and variance $\Sigma$ converges in distribution to a mean zero Gaussian random variable $W$ with variance $\Sigma$. In this setting, there do exist some results that bound $\Delta_f(W_n, W)$ for a finite sample size $n$. \citet[Chapter 13]{bhattacharya2010normal} is a classical reference for this when the random variables $X_i$ belong to the Euclidean space $\mathbb{R}^k$ and $f$ is an arbitrary Borel measurable function. Of course, $\Delta_f(W_n, W)$ cannot converge to zero for all Borel measurable functions and the bounds presented in Chapter 13 of~\cite{bhattacharya2010normal} involve an oscillation function of $f$ that implies a regularity condition on $f$. Formally, with $\|\cdot\|_2$ representing the Euclidean norm, define the oscillation function of $f:\mathbb{R}^d\to\mathbb{R}$ at a point $x$ with radius $\varepsilon$ as
\[
\omega_f(x; \varepsilon) := \sup\{|f(x) - f(y)|:\,\|x - y\|_2 \le \varepsilon\}. 
\]
One of the key quantities in the bound on $\Delta_f(W_n, W)$ is $\mathbb{E}[\omega_f(W, \varepsilon_n)]$ for some $\varepsilon_n$ converging to zero; see~\citet[Corollary 11.2, Thms 13.2 -- 13.3]{bhattacharya2010normal} for details. Unfortunately, the dependence on the dimension in these bounds is polynomial and hence, these results are not particularly useful for high-dimensional settings where the dimension can grow faster than the sample size. 

Moving beyond the Euclidean space ($\mathbb{R}^k$), for smooth (differentiable) functions $f$, one can bound $\Delta_f(W_n, W)$ even in Hilbert and Banach spaces. See, for example,~\cite{bentkus1993smoothness,rachev1989rates} and references therein. 

What is lacking in the literature, to our knowledge, is a bound on $\Delta_f$ that is locally adaptive. This means that the bound holds for all Borel measurable functions but at the same time yields the correct rate bound if the function $f$ is intrinsically low-dimensional or highly smooth. This is the main motivation for the current work. The idea we propose is that if the function $f$ can be written as $f(x) = \sum_{j=1}^{\infty} f_j(x)$ for a sequence of functions $f_j, j\ge1$ with $R_J = \sup_x|\sum_{j=J}^{\infty} f_j(x)| \to 0$ as $J\to\infty$, then we get $\Delta_f \le \sum_{j=1}^J \Delta_{f_j} + R_J$. By minimizing the bound over all $J$, one can bound $\Delta_f$. The resulting bound need not be optimal in terms of the dependence on $n$ with this approach, but they will often be adaptive. Moreover, most of the techniques discussed herein are applicable to non-iid random variables.

\paragraph{Organization.} The remaining article is organized as follows. In Section~\ref{sec:preliminaries}, we summarize the results of Part I along with some more exploration of those results. In Section~\ref{sec:integral-representation}, we describe and review integral representations for functions under several regularity conditions. In Section~\ref{sec:non-uniform-expansions}, we discuss non-uniform Edgeworth-type expansions for sums of independent univariate random variables and apply these results to the integral representations. In Section~\ref{sec:discussion}, we summarize the article and mention potential future directions.
\section{Preliminaries}\label{sec:preliminaries}
In part I~\citep{banerjee2023central},\footnote{From now,~\cite{banerjee2023central} is referred to as part I.} we have developed several simple results that can bound $\Delta_f$. In this section, we recall those results briefly. For any Borel measurable function $f:\mathcal{W}\to\mathbb{R}$, 
\begin{equation}\label{eq:expectation_and_level_sets}
\Delta_f(W_n, W) = \left|\int_{-\infty}^{\infty} [\mathbb{P}(W_{n} \in \mathcal{U}_{f,t}) - \mathbb{P}(W \in \mathcal{U}_{f, t})]\right|,
\end{equation}
where $\mathcal{U}_{f, t} = \{w\in\mathcal{W}:\, f(w) \ge t\}$ is the upper level set of $f$. Equality~\eqref{eq:expectation_and_level_sets} implies that one can control $\Delta_f(W_n, W)$ in terms of the difference between probabilities of the level sets. Differences between probabilities for sums of independent/dependent random variables/vectors can be controlled using traditional Berry--Esseen bounds. For example, the results of~\cite{bentkus2003dependence,bentkus2005lyapunov},~\cite{yaroslavtseva2008non},~\cite{raivc2019multivariate},~\cite{gotze1986rate}, and~\cite{paulauskas_Rackauskas_2012}, among others, provide such results for a large class of sets. (The last two references also contain results specific to lower-level sets of smooth functions.)

It is not difficult to find functions whose level sets belong to these favorable classes. Consider the following simple examples.
\paragraph{Convex sets.} If $f:\mathcal{W}\to\mathbb{R}$ is a quasiconcave function (i.e., $f(\lambda w_1 + (1 - \lambda)w_2) \ge f(w_1)\wedge f(w_2)$), then $\mathcal{U}_{f, t}$ is a convex set. In fact, this is the characterization of quasiconcave functions~\citep{diewert1981nine}. \cite{michel1978remark} provides bounds on $\mathbb{P}(W_n\in A) - \mathbb{P}(W\in A)$ that depend on $A$ and converge to zero as $n\to\infty$, if $W_n$ is an average of $n$ independent random vectors in $\mathbb{R}^d$ and $W$ is a Gaussian random vector $(W\sim N(\mathbb{E}[W_n], \mbox{Var}(W_n)))$. A similar result is provided in~\cite{hipp1979convergence} if $W_n$ is an average of mixing random vectors. Also, see~\cite{rotar1970non},~\cite{sazonov1974new},~\cite{senatov1982some},~\cite{fomin1983estimate}, and~\cite{jirak2015non} for more results of this flavor.
\paragraph{Euclidean balls.} If $f(w) = g(\|w - a\|_2), w\in\mathbb{R}^d$ for some non-increasing function $g$, then the upper level sets of $f$ are Euclidean balls centered at $a\in\mathbb{R}^d$\footnote{$\|\cdot\|_2$ throughout represents the Euclidean norm}. \cite{bogatyrev2006non} provides bounds on $\mathbb{P}(W_n \in A) - \mathbb{P}(W\in A)$ that depend on Euclidean ball $A$ and converge to zero as $n\to\infty$, if $W_n$ is an average of $n$ independent random variables in a Hilbert space and $W$ is the corresponding Gaussian random variables; here $\|\cdot\|_2$ represents the Hilbert space norm. Further results can be found in~\cite{ulyanov1986normal},~\cite{yurinskii1983accuracy},~\cite{tikhomirov1994accuracy},~\cite{senatov2011normal} and~\cite{paulauskas_Rackauskas_2012}. The class of all functions with Euclidean balls centered at $a$ as upper level sets is the same as the class of all functions of the form $w\mapsto g(\|w - a\|_2)$ for some non-increasing function $g$. 
\paragraph{Half-spaces.} If $f(w) = g(\langle w, a\rangle - b), w\in\mathcal{W}$ for some inner product space $\mathcal{W}$ and some non-increasing function $g$, then the upper level sets of $f$ are half-spaces (i.e., of the form $\mathcal{U}_{f,t} := \{w\in\mathbb{R}^d:\, \langle a, w\rangle \ge b + g^{-1}(t)\}$). Bounds on the difference of probabilities for half-spaces can be obtained from univariate Berry--Esseen bounds even if $W_n, W$ are multivariate/Hilbert/Banach space valued random variables. This is because $\mathbb{P}(W_n \in \mathcal{U}_{f, t}) = \mathbb{P}(\langle a, W_n\rangle \ge b + g^{-1}(t))$ and hence, if $W_n$ is an average of random variables in some measurable space, then $\langle a, W_n\rangle$ is an average of univariate random variables depending on $a$; although a simple fact, this was also mentioned in the remark following Corollary 3 of~\cite{paulauskas1976convergence}. Such bounds for independent or dependent data can be found in the literature; see, for example,~\cite{michel1978remark},~\cite{heinrich1985non}, and~\cite{jirak2015non}. In all these cases, one can also obtain asymptotic expansions to get precise bounds on $\Delta_f$. We will present an analysis of this kind in the following sections.

The discussion above allows us to get bounds for functions whose level sets belong to a favorable class. This may not always be the case. For example, sum of two quasiconcave functions need not be quasiconcave and hence, the upper level sets of the sum of two quasiconcave functions need not be convex. However, we can use the simple fact that $\Delta_{f_1 + f_2}(W_n, W) \le \Delta_{f_1}(W_n, W) + \Delta_{f_2}(W_n, W)$ to get a bound. This relation need not be restricted to the sum of two functions but can be extended to the sum of uncountably many functions. Formally, we have the following result (proved in Appendix~\ref{sec:proof-signed-measure-bound}).
\begin{proposition}\label{prop:signed-measure-bound}
    Suppose $\{f_{\lambda}:\,\lambda\in\Lambda\}$ is a parametrized class of functions and $\mu(\cdot)$ is a finite signed measure on $\Lambda$. Then, for $g(x) = \int_{\Lambda} f_{\lambda}(x)\mu(d\lambda)$ and for any random variables $W_n, W$, we have
    \[
    \Delta_g(W_n, W) ~\le~ \int_{\Lambda} \left|\Delta_{f_{\lambda}}(W_n, W)\right||\mu|(d\lambda),
    \]
    where $|\mu|(\cdot)$ is the variation of $\mu$.
\end{proposition}
There are several simple applications of Proposition~\ref{prop:signed-measure-bound}, even with discrete measures $\mu$ supported on finite or countable sets. For example, it is well-known that any differentiable function with an $L$-Lipschitz derivative can be written as a difference between a convex function with $2L$-Lipschitz derivative and a quadratic function~\citep{zlobec2006characterization}.

Proposition~\ref{prop:signed-measure-bound}, in particular, applies to functions of the form $f(x) = \sum_{j=1}^{\infty} \theta_j\phi_j(x)$ for some coefficients such that $\sum_{j=1}^{\infty} |\theta_j| < \infty$. This is the classical case of basis expansion and for functions defined on Euclidean spaces, it is well-known that Fourier/spline/wavelet bases can approximate any function in $L_2$, at least for functions with bounded support. In these cases, however, the coefficients are only implied to be square summable and not absolutely summable. For Haar bases (which is an example of wavelet bases), $\phi_j(x)$ has level sets that are hyperrectangles, and for hyperrectangles, several bounds for the difference of probabilities exist~\citep{chernozhukov2023nearly,bong2022high,fang2023high}.  

Proposition~\ref{prop:signed-measure-bound} plays a crucial role in this paper and we will present several applications of it where $f_{\lambda}(\cdot)$ are functions whose level sets are either Euclidean balls or half-spaces. In the following section, we provide examples of integral representations of the type $g(x) = \int_{\Lambda} f_{\lambda}(x)\mu(d\lambda)$ based on the Fourier transform of $g$.
\section{Integral Representations}\label{sec:integral-representation}
In this section, we present some sufficient conditions under which functions can be represented as finite integrals of functions whose level sets belong to a favorable class (e.g., Euclidean balls or half-spaces).
\subsection{Norm balls}
In this section, we use the approximation of functions $f$ by integrals of radial functions (i.e., functions whose level sets are Euclidean balls) to bound $\Delta_f(W_n, W)$. It is a well-known fact that any function can be approximated by convolving the function with a mollifier (i.e., an infinitely differentiable function that approximates the Dirac delta). Formally, in the finite-dimensional Euclidean space $\mathbb{R}^d$, any function $f:\mathbb{R}^d\to\mathbb{R}$ can be approximated by $f_h$ where $f_h(x) = h^{-d}\int_{\mathbb{R}^d} K(x/h, y/h)f(y)dy$, for any function $K(\cdot, \cdot)$. In particular, Proposition~4.3.31 of~\cite{gine2021mathematical} states that if $\int_{\mathbb{R}^d}\sup_{v\in\mathbb{R}^d}|K(v, v-u)|du < \infty$ and $\int_{\mathbb{R}^d} K(x, y)dy = 1$ for all $x\in\mathbb{R}^d$, then $\|f_h - f\|_{\infty} = \sup_{x\in\mathbb{R}^d}|f_h(x) - f(x)|$ converges to zero. We can use this approximation result with Proposition~\ref{prop:signed-measure-bound} to get a bound for $\Delta_f(W_n, W)$. Note that
\[
\left|\mathbb{E}[f(W_n)] - \mathbb{E}[f_h(W_n)]\right| \le \mathbb{E}|f(W_n) - f_h(W_n)| \le \|f_h - f\|_{\infty},  
\]
and
\[
\mathbb{E}[f_h(W_n)] = \frac{1}{h^d}\int_{\mathbb{R}^d} \mathbb{E}\left[K\left(\frac{W_n}{h}, \frac{y}{h}\right)\right]f(y)dy,
\]
assuming that the right hand side exists.
Therefore, 
\begin{equation}\label{eq:Delta-f-inequality}
\Delta_f(W_n, W) \le 2\|f_h - f\|_{\infty} + \Delta_{f_h}(W_n, W).
\end{equation}
Because the right hand side depends on the kernel $K(\cdot, \cdot)$ and bandwidth $h$, but the left hand side does not, we get 
\begin{equation}\label{eq:Delta-f-inequality-with-infimum}
\Delta_f(W_n, W) \le \inf_{K, h}\,\left\{2\|f_h - f\|_{\infty} + \Delta_{f_h}(W_n, W)\right\}.
\end{equation}
Proposition 4.3.33 of~\cite{gine2021mathematical} also states that if $f\in C^m(\mathbb{R}^d)$ (the space of all functions whose $\lfloor m\rfloor$-th derivative is bounded and is $m-\lfloor m\rfloor$-Holder continuous), then there exists (higher-order) kernel such that $\|f_h - f\|_{\infty} \le C_{f}h^m$; a higher-order kernel with order $\ell$ means that $\int_{\mathbb{R}^d} K(v, u + v)u^{\alpha}du = 0$ for all $v\in\mathbb{R}^d$ for every monomial $u^{\alpha}$ with total degree $|\alpha|$ less than $\ell$. Such higher-order kernels can be constructed easily using a technique called twicing~\citep{stuetzle2006some}: if $K$ is a kernel of order $\ell$, then $K_2 = 2K - K*K$ is a kernel of order $2\ell$.\footnote{We use the notation '$*$' to denote convolution.} 

Calculating the infimum in~\eqref{eq:Delta-f-inequality-with-infimum} can be difficult in general. To illustrate the usefulness, we consider the Gaussian kernel to obtain a concrete bound, i.e., $K(u, v) = (2\pi)^{-d/2}\exp(-\|u - v\|_2^2/2)$. With this Gaussian kernel, we have
\[
\mathbb{E}[f_h(W_n)] = \frac{1}{h^d(2\pi)^{d/2}}\int_{\mathbb{R}^d} \mathbb{E}\left[\exp\left(-\frac{\|W_n - y\|_2^2}{2h}\right)\right]f(y)dy.
\]
Applying~\eqref{eq:expectation_and_level_sets}, we get
\[
\mathbb{E}\left[\exp\left(-\frac{\|W_n - y\|_2^2}{h}\right)\right] = \int_0^1 \mathbb{P}\left(\|W_n - y\|_2 \le \sqrt{2h\log(1/t)}\right)dt.
\]
Therefore, 
\[
\Delta_{f_h}(W_n, W) \le \frac{(2\pi)^{-d/2}}{h^d}\int_{\mathbb{R}^d\times[0, 1]} \left|\mathbb{P}\left(\frac{\|W_n - y\|_2}{\sqrt{2h\log(1/t)}} \le 1\right) - \mathbb{P}\left(\frac{\|W - y\|_2}{\sqrt{2h\log(1/t)}} \le 1\right)\right|f(y)dtdy. 
\]
This inequality allows us to control $\Delta_{f_h}(W_n, W)$ using Berry--Esseen bounds for Euclidean balls. The calculations presented in this section can be easily generalized to Hilbert space random variables. As an example, we now proceed to control the right hand side when $W_n$ is a scaled average of $n$ iid random vectors in $\mathbb{R}^d$ for $d \ge 6$. Following the results of~\cite{senatov1992some,senatov1993estimates}, we obtain that if $W_n = n^{-1/2}\sum_{i=1}^n X_i$ for mean zero random vectors $X_i\in\mathbb{R}^d$, then with $r_h(y) = |\sqrt{2h\log(1/t)} - \|y\|_2|$,
\begin{align*}
&\left|\mathbb{P}\left(\frac{\|W_n - y\|_2}{\sqrt{2h\log(1/t)}} \le 1\right) - \mathbb{P}\left(\frac{\|W - y\|_2}{\sqrt{2h\log(1/t)}} \le 1\right)\right|\\ 
&\quad\le \frac{C\beta_3/n^{1/2}}{1 + r_h^3(y)/\sigma^3}\left\{\frac{(2h\log(1/t))^{3/2}}{\sigma_1\cdots\sigma_6}\exp\left(-\frac{cr_{h}^2(y)}{\sigma^2}\right) + \frac{1}{\sigma^3} + \frac{1}{\sqrt{\sigma_1\cdots\sigma_6}}\exp\left(-\frac{cr_h^2(y)}{\sigma^2}\right)\right\},
\end{align*}
where $\beta_3 = \mathbb{E}[\|X\|_2^3]$, $\sigma^2 = \mathbb{E}[\|X\|_2^2]$, and $\sigma_1^2 \ge \sigma_2^2 \ge \cdots \ge \sigma_6^2 > 0$ are the largest $6$ eigenvalues of $\Sigma = \mathbb{E}[XX^{\top}].$
Clearly, if $\|y\|_2 \ge 2\sqrt{2h\log(1/t)}$, then $r_h(y) \ge \|y\|_2/2$ and hence, the difference between the probabilities can be bounded by
\[
\frac{C\beta_3n^{-1/2}}{(1 + \|y\|_2^3\mathbf{1}\{\|y\|_2 \ge \sqrt{8h\log(1/t)}\}/(8\sigma^3))}\left\{\frac{(2h\log(1/t))^{3/2}}{(\sigma_1\cdots\sigma_6)} + \frac{1}{\sigma^{3}} + \frac{1}{(\sigma_1\cdots\sigma_6)^{1/2}}\right\}.
\]
Therefore, to bound $\Delta_{f_h}(W_n, W)$, it suffices to compute the integral of the right hand side weighted by $f(y)$ over $\mathbb{R}^d\times[0, 1].$ A simple calculation (presented in Appendix~\ref{sec:proof-of-norm-ball-result}) yields
\begin{align*}
&\Delta_{f_h}(W_n, W)\\ 
&\le \frac{C\beta_3}{n^{1/2}(2\pi)^{d/2}}\left[\frac{h^{3/2}}{(\sigma_1\cdots\sigma_6)} + \frac{1}{\sigma^{3}} + \frac{1}{(\sigma_1\cdots\sigma_6)^{1/2}}\right]\left(\int_{\mathbb{R}^d} \left\{\frac{8\sigma^3}{8\sigma^3 + h^3\|y\|_2^3}+\exp\left(-\frac{\|y\|_2^2}{16}\right)\right\}f(yh)dy\right).
\end{align*}
As long as the integrals on the right hand side are finite (which holds, for example, if $f$ belongs to Schwartz class), the right hand side converges to zero as $n\to\infty$ for every $h > 0$. In particular, if we have $g(x) = \int \exp(\|x-y\|_2^2/2h)f(y)dy/(\sqrt{2\pi}h)^d$, then $\Delta_g(W_n, W) = \Delta_{f_h}(W_n, W)$ and so the bound above applies. The same bound also applies to higher-order twicing kernels starting with a Gaussian kernel, because the convolution of a Gaussian with itself is another Gaussian. Finally, one can extend the result above to obtain bounds for $\Delta_g(W_n, W)$ for functions $g$ of the form $g(x) = \int_{\mathbb{R}^d}\int_0^{\infty} h^{-d}K(\|x - y\|_{2}/h)f(y)\pi(h)dydh$, for some function $\pi(\cdot)$; this follows by integrating the above bound with respect to $\pi(h)$. For a class of functions that satisfy such a representation for some choices of $K(\cdot)$ and $\pi(\cdot)$, see~\cite{girosi1992convergence} and~\cite{girosi1994regularization}. In general, for any function $f$, we apply inequality~\eqref{eq:Delta-f-inequality} and minimize over $h > 0$.

We end this section with a discussion of the relation to reproducing kernel Hilbert spaces (RKHS). For any positive definite function $K:\mathcal{X}\times\mathcal{X}\to\mathbb{R}$~\citep[Definition 12.6]{wainwright2019high}, the RKHS corresponding to $K$ is defined as the unique Hilbert space $\mathbb{H}_K = \{f:\, f(x) = \int_{\mathcal{X}} K(x, y)f(y)dy\mbox{ }\forall\mbox{ }x\in\mathcal{X}\}$. Here $\mathcal{X}$ is an arbitrary metric space. For any function $f\in\mathbb{H}_K$, the techniques developed in this section can be readily applied. In particular, we have
\begin{equation}\label{eq:RKHS-bound}
\begin{split}
\Delta_f(W_n, W) &= \left|\int_{\mathcal{X}} \mathbb{E}[K(W_n, y) - K(W, y)]f(y)dy\right|,\\
&\le \int_{\mathcal{X}} \left|\mathbb{E}[K(W_n, y) - K(W, y)]\right||f(y)|dy.
\end{split}
\end{equation}
The RKHS induced by the Gaussian kernel $K(x, y) = \exp(-\|x - y\|_2^2/\gamma)$ is explicitly described in~\citet[Theorem 1]{minh2010some}. For all such functions, the bound $\Delta_{f_h}(W_n, W)$ can be used with a fixed $h$; it is also known that the Gaussian RKHS is dense in $L_p$ for any $p\ge1$; see~\citet[Sec. 4.6 and Thm. 4.63]{steinwart2008support}. Furthermore, the classical smoothness class of Sobolev space $W_2^m(\mathbb{R}^d)$ is also an RKHS with $m > d/2$~\citep{novak2018reproducing}. The corresponding kernel can also be written as a radial function that involves Bessel functions of the third kind; see~\citet[Theorem 9.9]{schaback2007kernel} for more details. Also, see~\cite{de2021reproducing} for extensions of RKHS defined on manifolds. 

What is the main takeaway from the results in this section? We have provided a general recipe to obtain bounds for $\Delta_f(W_n, W)$ for general $f$. If $f$ is approximable by a smooth function, then one can get good bounds on $\Delta_f(W_n, W)$. In particular, when $f$ belongs to an RKHS with radial kernel $K(\cdot, \cdot)$, then the bounds obtained on $\Delta_f(W_n, W)$ do not have excessive dependence on the dimension and on the covariance matrix of the underlying random variables. In fact, the explicit dependence on the dimension $d$ actually decreases with increasing dimension. The results also apply to functions defined on a Hilbert space.
\subsection{Half-spaces}
In this section, we provide integral representations of functions $f$ in terms of functions whose level sets are half-spaces. As mentioned before, in this case, we can obtain results for high-/infinite-dimensional random variables using results for univariate random variables. Following part I, we start with results from neural networks literature to provide sufficient conditions on functions for integral representations.

We start with the recent work of~\cite{klusowski2018approximation} and extend it to present integral representations for functions defined on $\mathbb{R}^d$ (instead of $[-1, 1]^d$ used in~\cite{klusowski2018approximation}). In this section, we restrict to functions $f$ for which there is a Fourier representation of the form 
\[
f(x) = \int_{\mathbb{R}^d} e^{i\langle \omega, x\rangle}\widecheck{f}(\omega)d\omega,
\]
for some complex function $\widecheck{f}(\cdot)$. In what follows, we show that under certain weighted integrability assumptions on $\widecheck{f}$, $f$ satisfies an integral representation in terms of functions whose level sets are half-spaces. Theorem 2 of~\cite{klusowski2018approximation} proved that if $f$ is defined on $[-1, 1]^d$ and $v_{f,2} = \int_{\mathbb{R}^d}\|\omega\|_1^2|\widecheck{f}(\omega)|d\omega < \infty$, then there exists a probability measure $P$ on $\{-1, 1\}\times[0, 1]\times\mathbb{R}^d$,
\[
f(x) = f(0) + \langle x, \nabla f(0)\rangle - v\int_{\{-1, 1\}\times [0, 1]\times \mathbb{R}^d} (z\langle a, x\rangle - t)_+\mbox{sgn}(\cos(\|\omega\|_1zt + b(\omega)))dP(z, t, \omega),
\]
for all $x\in[-1, 1]^d$, where $b(\omega)$ is defined by the relation $\widecheck{f}(\omega) = e^{ib(\omega)}|\widecheck{f}(\omega)|$. The choice $\|\cdot\|_1$-norm here is rather arbitrary for the purpose of central limit theorems (but is well-motivated in the context of neural network learning). The predecessor work~\cite{barron1993universal} actually uses $\|\cdot\|_2$-norm in place of the $\|\cdot\|_1$-norm and derives a similar representation when $v_{f, 1}^{(2)} = \int_{\mathbb{R}^d}\|\omega\|_2|\widecheck{f}(\omega)|d\omega < \infty$. The following theorem presents a minor extension of the results of~\cite{klusowski2018approximation} that does not require the specification of the norm or restrict the domain.
\begin{theorem}
\label{thm:integral-representation-ridge}
Let $f:\mathbb{R}^d\to\mathbb{R}$ be a function with Fourier transform $\widecheck{f}:\mathbb{R}^d\to\mathbb{R}$. Then for all $x\in\mathbb{R}^d$ such that 
\begin{equation}\label{eq:ReLU-representation-assumption}
\int_{\mathbb{R}^d} \min\left\{2|\langle \omega, x\rangle|,\,|\langle \omega, x\rangle|^2/2\right\}|\widecheck{f}(\omega)|d\omega ~<~ \infty,
\end{equation}
we have
\[
f(x) = f(0) + f'(0)[x] - \int_{\mathbb{R}^d}\int_0^{\infty} \mathbb{E}_{\varepsilon}[(\varepsilon\langle\omega,x\rangle - u)_+\cos(u + \varepsilon b(\omega))]|\widecheck{f}(\omega)|dud\omega,
\]
where $b(\cdot)$ is defined via $\widecheck{f}(\omega) = |\widecheck{f}(\omega)|e^{ib(\omega)}$ and $\mathbb{E}_{\varepsilon}[g(\varepsilon, u, \omega)] = (g(1, u, \omega) + g(-1, u, \omega))/2$ for any function $g$. Moreover, if for any two random variables $W, W_n\in\mathbb{R}^d$ with $\mathbb{E}[W] = \mathbb{E}[W_n] = 0$,
\begin{equation}\label{eq:random-variable-ReLU-representation-assumption}
\begin{split}
\int_{\mathbb{R}^d} \mathbb{E}[\min\{2|\langle\omega, W_n\rangle|, |\langle\omega, W_n\rangle|^2/2\}]|\widecheck{f}(\omega)|d\omega &< \infty,\\
\int_{\mathbb{R}^d} \mathbb{E}[\min\{2|\langle\omega, W\rangle|, |\langle\omega, W\rangle|^2/2\}]|\widecheck{f}(\omega)|d\omega &< \infty,
\end{split}
\end{equation}
then
\begin{equation}\label{eq:ReLU-bound}
\begin{split}
\Delta_f(W, W_n) &= \left|\int_{\mathbb{R}^d}\int_0^{\infty} \mathbb{E}_{\varepsilon}[\mathbb{E}\{(\varepsilon\langle\omega, W_n\rangle - u)_+ - (\varepsilon\langle\omega, W\rangle - u)_+\}\cos(u + \varepsilon b(\omega))]|\widecheck{f}(\omega)|dud\omega\right|,\\
&\le \mathbb{E}_{\varepsilon}\left[\int_{\mathbb{R}^d}\int_0^{\infty} |\mathbb{E}[(\varepsilon\langle\omega, W_n\rangle - u)_+ - (\varepsilon\langle\omega,W\rangle - u)_+]||\widecheck{f}(\omega)|dud\omega\right].
\end{split}
\end{equation}
\end{theorem}
Theorem~\ref{thm:integral-representation-ridge} is stated to represent (or bound) $\Delta_f(W_n, W)$ in terms of $\Delta_g(W_n, W)$ where $g(x) = \langle \omega, x\rangle$ for some $\omega\in\mathbb{R}^d$. Because $g(W_n)$ is an average of univariate random variables if $W_n$ is an average of random variables in $\mathbb{R}^d$, we can control $\Delta_g(W_n, W)$ easily even for dependent random variables, that too without any explicit dimension dependence. The condition~\eqref{eq:ReLU-representation-assumption} stems from the inequality $|e^{-iz} - iz - 1| \le \min\{|z|, |z|^2/2\}.$ Note that assumption~\eqref{eq:ReLU-representation-assumption} implies that the function $f$ is differentiable at zero. One can also obtain a result similar to Theorem~\ref{thm:integral-representation-ridge} with $(\varepsilon\langle\omega,x\rangle - u)_+$ replaced by $(\varepsilon\langle\omega,x\rangle - u)_+^s$ for $s \ge 0$ by replacing the assumption~\eqref{eq:ReLU-representation-assumption} with $\int_{\mathbb{R}^d} \min\{|\langle\omega, x\rangle|^s, |\langle\omega, x\rangle|^{s+1}\}|\widecheck{f}(\omega)|d\omega < \infty$. This revised assumption implies that the function is $s$-times differentiable at zero.

In what follows, we provide further integral representations for functions under regularity conditions akin to~\eqref{eq:ReLU-representation-assumption}. The following result is modeled after the result of~\cite{irie1988capabilities}; also see Eqs. (19)-(21) of~\cite{siegel2020approximation}.
\begin{theorem}\label{thm:integral-representation-general-activation}
Suppose $h:\mathbb{R}\to\mathbb{R}$ is an integrable function with Fourier transform $\widecheck{h}(\cdot)$ such that $\widecheck{h}(a) \neq 0$ for some $a\in\mathbb{R}$. Then any function $f:\mathbb{R}^d\to\mathbb{R}$ with an integrable Fourier transform $\widecheck{f}(\cdot)$ can be represented as
\[
f(x) = \int_{\mathbb{R}^{d}}\int_{\mathbb{R}}\frac{1}{2\pi|\widecheck{h}(a)|}h\left(a^{-1}\langle\omega,x\rangle + u\right)|\widecheck{f}(\omega)|\cos(au + c_h - b(\omega))dbd\omega,\quad\mbox{for all}\quad x\in\mathbb{R}^d,
\]
where $\widecheck{h}(a) = |\widecheck{h}(a)|e^{ic_h}$, and $\widecheck{f}(\omega) = |\widecheck{f}(\omega)|e^{ib(\omega)}$.
\end{theorem}
Theorem~\ref{thm:integral-representation-general-activation} can be used in the same way as Theorem~\ref{thm:integral-representation-ridge} to get bounds on $\Delta_f$. In this result also, we only need to study $\mathbb{E}[h(\langle \omega, W_n\rangle + u) - h(\langle \omega, W\rangle + u)]$ for fixed $\omega$ and $u$. If $W_n$ is an average of random variables in $\mathbb{R}^d$, then $\langle\omega,W_n\rangle$ is also an average of univariate random variables and hence, univariate central limit theorems can be used to bound $\Delta_f$ for arbitrary dimension $d\ge1$. A limitation of Theorem~\ref{thm:integral-representation-general-activation} is that the activation function $h(\cdot)$ is required to be integrable and several commonly used functions such as logistic, ReLU, or Heaviside functions do not satisfy this condition. Interestingly, there is a simple way to rectify this limitation. Lemma 1 of~\cite{funahashi1989approximate} shows that for any non-constant, bounded, monotone increasing continuous function $\phi:\mathbb{R}\to\mathbb{R}$, $h(t) = \phi(t + \alpha) - \phi(t - \alpha)$ is an integrable function for any $\alpha > 0.$ Furthermore, there exists an $a \in\mathbb{R}$ such that $\widecheck{h}(a) \neq 0$. The advantage of using a monotonically increasing bounded activation function $h$ is that the right side of~\eqref{eq:expectation_and_level_sets} is a finite integral. A similar integral representation is also obtained in~\cite{makovoz1998uniform} under the assumption that $\sup_{u:\|u\|_2 = 1}\int_0^{\infty} r^{d}|\widecheck{f}(ru)|dr < \infty$; see Remark 1 of~\cite{klusowski2018approximation}.

Similar to the use of the Fourier transform, one can obtain integral representations of functions using other transforms such as the Radon transform. We present one such example below and leave a detailed study for future work. The following result is taken from~\cite{kainen2010integral}. For any $u\in\mathbb{R}^d$ such that $\|u\|_2 = 1$ and $b \in\mathbb{R}$, define
\[
H_{u,b}^{-} = \{x\in\mathbb{R}^d:\,\langle u,x\rangle + b \le 0\}.
\]
The linear operator Laplacian $\Delta$ is defined by $\Delta g = \sum_{j=1}^d \partial^2g/\partial x_j^2$, for a twice differentiable function $g:\mathbb{R}^d\to\mathbb{R}$. For a positive integer $m\ge1$, $\Delta^m$ denotes the Laplacian iterated $m$ times, while $\Delta^0$ is the identity operator. Define $k_d = 2\lceil(d+1)/2\rceil$. Call a function $f:\mathbb{R}^d\to\mathbb{R}$ to have controlled decay if $f$ is $k_d$-times continuously differentiable and for each multi-index $\alpha = (\alpha_1, \ldots, \alpha_d)$ with $\sum_{j=1}^d |\alpha_j| \le k_d$, there exists $\varepsilon > 0$, such that $\lim_{\|x\|_2\to\infty} \partial^{\alpha}f(x)\|x\|_2^{|\alpha| + \varepsilon} = 0$; this condition means that the derivatives of the order less than $k_d$ converge to zero at infinity at least as fast as a polynomial. For any function $f$ of controlled decay, define
\[
w_{f}(u,b) = a_d\times
\begin{cases}
\int_{H_{u,b}^{-}} \Delta^{k_d/2}f(y)dy,&\mbox{ if }d\mbox{ is odd,}\\
\int_{\mathbb{R}^d} \Delta^{k_d/2}f(y)\alpha(\langle u, y\rangle + b)dy, &\mbox{if }d\mbox{ is even,}
\end{cases}
\]
where $a_d = (-1)^{(d-1)/2}/(2(2\pi)^{d-1})$ if $d$ is odd and $a_d = (-1)^{(d-2)/2}/(2\pi)^d$ if $d$ is even, and $\alpha(t) = t\log(e/|t|)$ for $t\neq 0$ with $\alpha(0) = 0$. With these notations at hand, Theorem~4.2 of~\cite{kainen2010integral} states that any function $f$ of controlled decay satisfies
\begin{equation}\label{eq:integral-representation-Radon-transform}
f(x) = \int_{S^{d-1}\times\mathbb{R}} w_f(u, b)\mathbf{1}\{\langle u, x\rangle + b \ge 0\}dudb,
\end{equation}
where $S^{d-1} = \{u\in\mathbb{R}^d:\,\|u\|_2 = 1\}$. For more representations related to Radon transform, see~\cite{petrosyan2020neural} and~\cite{abdeljawad2022integral}. 

Inequalities of the type~\eqref{eq:RKHS-bound} and~\eqref{eq:ReLU-bound} allow one to easily bound $\Delta_f(W_n, W)$ for a large class of functions 
(in fact, many of which are dense in $L_2$) easily in terms of $\Delta_g(W_n, W)$ for a small class of functions $g$. 
(For~\eqref{eq:RKHS-bound}, it suffices to bound $\mathbb{E}[K(W_n, y) - K(W, y)]$ for each $y$ and for~\eqref{eq:ReLU-bound}, it suffices to bound 
$\mathbb{E}[(\varepsilon\langle\omega, W_n\rangle - u)_+ - (\varepsilon\langle\omega, W\rangle - u)_+]$ for each $\varepsilon\in\{-1, 1\}, u\in[0, \infty)$, and $\omega\in\mathbb{R}^d$). Moreover, the equalities in~\eqref{eq:RKHS-bound},~\eqref{eq:ReLU-bound}, Theorem~\ref{thm:integral-representation-general-activation}, and~\eqref{eq:integral-representation-Radon-transform} also allow one to get precise asymptotics of $\Delta_f(W_n, W)$ via upper and lower bounds.
\section{Expansions and Applications}\label{sec:non-uniform-expansions}
In Section~\ref{sec:integral-representation}, we provided several integral representations for functions in terms of Euclidean balls and half-spaces. In this section, we review non-uniform expansions for the difference of probabilities for averages of univariate random variables and consider the application of these expansions for precise bounds on $\Delta_f$. In the discussion to follow, we restrict ourselves to independent and identically distributed random variables. Extensions to non-independent and/or non-identically distributed random variables are definitely of interest but will be explored elsewhere.

Let $W_1, \ldots, W_n$ be a sequence of independent and identically distributed real-valued random variables with mean zero and variance $\sigma^2 > 0$. Set $v(b) = \mathbb{E}[e^{ibW_1}]$, and $\beta_3 = \mathbb{E}[|W_1|^3]$. Theorem 5.18 of~\cite{Petrov1995} (with $k = 3$) implies the following result.
\begin{theorem}\label{thm:Edgeworth-non-uniform-expansion}
Let $Z\sim N(0, 1)$. If $\beta_3 < \infty$, then there exists a constant $C > 0$ such that for all $x\in\mathbb{R}$,
\begin{equation}\label{eq:edgeworth-non-uniform}
\begin{split}
&\left|\mathbb{P}\left(\frac{1}{\sqrt{n\sigma^2}}\sum_{i=1}^n W_i \le x\right) - \mathbb{P}(Z \le x) - \frac{(1 - x^2)e^{-x^2/2}\mathbb{E}[W_1^3]}{6\sqrt{2\pi n}\sigma^3}\right|\\ 
&\quad\le C\frac{\mathbb{E}[|W_1|^3\mathbf{1}\{|W_1| \ge \sigma(1 + |x|)n^{1/2}\}]}{\sigma^3n^{1/2}(1 + |x|)^3}\\
&\qquad+ C\frac{\mathbb{E}[|W_1|^4\mathbf{1}\{|W_1| < \sigma(1 + |x|)n^{1/2}\}]}{\sigma^4n(1 + |x|)^{4}}\\
&\qquad+ C\left(\sup_{|b| \ge \sigma^2/(12\beta_3)}|v(b)| + \frac{1}{2n}\right)^n\frac{n^{6}}{(1 + |x|)^4}.
\end{split}
\end{equation}
\end{theorem}
The assumption of $\beta_3 < \infty$ implies that the right hand side converges to zero as $n\to\infty$ as long as $\sup_{|b| \ge \sigma^2/(12\beta_3)}|v(b)| < 1$, but cannot guarantee any specific rate of convergence. Stronger assumptions such as $\mathbb{E}[|W_1|^{3 + \gamma}] < \infty$ can allow one to obtain a specific decay rate (with respect to $n$) for the right hand side. In particular, if $\mathbb{E}[|W_1|^4] < \infty$, then
\begin{equation}\label{eq:Edgeworth-non-uniform-fourth-moment}
\begin{split}
&\left|\mathbb{P}\left(\frac{1}{\sqrt{n\sigma^2}}\sum_{i=1}^n W_i \le x\right) - \mathbb{P}(Z \le x) - \frac{(1 - x^2)e^{-x^2/2}\mathbb{E}[W_1^3]}{6\sqrt{2\pi n}\sigma^3}\right|\\ 
&\quad\le C\frac{\mathbb{E}[|W_1|^4]}{\sigma^4n(1 + |x|)^4} + C\left(\sup_{|b| \ge \sigma^2/(12\beta_3)}|v(b)| + \frac{1}{2n}\right)^n\frac{n^{6}}{(1 + |x|)^4}.
\end{split}
\end{equation}
The advantage of the results of the type~\eqref{eq:edgeworth-non-uniform} or~\eqref{eq:Edgeworth-non-uniform-fourth-moment} is that they provide a precise description of the difference of probabilities along with a bound that depends on $x$. These features will turn out to be very useful when applying them to the integral representations from the previous section.
\subsection{Application to ReLU functions}\label{sec:ReLU-Edgeworth}
For any random variable $U$, we have
\[
\mathbb{E}[(U - t)_+] = \int_0^{\infty} \mathbb{P}((U - t)_+ > s) = \int_0^{\infty} \mathbb{P}(U > t + s)ds = \int_0^{\infty} (1 - F_U(s + t))ds,
\]
where $F_U(\cdot)$ is the cumulative distribution function of $U$. Therefore,
\begin{align*}
\Delta_{n,\mathrm{ReLU}}(t) &:= \mathbb{E}\left(\frac{1}{\sqrt{n\sigma^2}}\sum_{i=1}^n W_i - t\right)_+ - \mathbb{E}[(Z - t)_+]\\
&= \int_0^{\infty} \left[\mathbb{P}(Z \le t + s) - \mathbb{P}\left(\frac{1}{\sqrt{n\sigma^2}}\sum_{i=1}^n W_i \le t + s\right)\right]ds.
\end{align*}
Inequality~\eqref{eq:Edgeworth-non-uniform-fourth-moment} now implies that
\begin{align*}
&\left|\Delta_{n,\mathrm{ReLU}}(t) - \frac{\mathbb{E}[W_1^3]}{6\sqrt{2\pi n}\sigma^3}\int_0^{\infty} {(1 - (t+s)^2)e^{-(t + s)^2/2}}ds\right|\\ 
&\quad\le C\left[\frac{\mathbb{E}[|W_1|^4]}{\sigma^4n} + n^6\left(\sup_{|b| \ge \sigma^2/(12\beta_3)}|v(b)| + \frac{1}{2n}\right)^n\right]\int_0^{\infty} \frac{1}{(1 + |t+s|)^4}ds.
\end{align*}
A similar result can be obtained from Theorem~\ref{thm:Edgeworth-non-uniform-expansion}. Simplifying the bound above, we have proved the following result.
\begin{proposition}\label{prop:Delta-ReLU-bound}
Under the notation above, we have
\begin{equation}\label{eq:Delta-ReLU-each-t}
\left|\Delta_{n,\mathrm{ReLU}}(t) + \frac{te^{-t^2/2}\mathbb{E}[W_1^3]}{6\sqrt{2\pi n}\sigma^3}\right| \le C\left[\frac{\mathbb{E}[|W_1|^4]}{\sigma^4n} + n^6\left(\sup_{|b| \ge \sigma^2/(12\beta_3)}|v(b)| + \frac{1}{2n}\right)^n\right]\kappa(t),
\end{equation}
where
\[
\kappa(t) := \begin{cases}
(3(1 + t)^3)^{-1}, &\mbox{if }t \ge 0,\\
2/3 - (3(1-t)^3)^{-1}, &\mbox{if }t < 0.
\end{cases}.
\]
Moreover, 
\begin{equation}\label{eq:Delta-ReLU-integral}
\begin{split}
&\left|\int_0^{\infty} \Delta_{n,\mathrm{ReLU}}(t)dt + \frac{\mathbb{E}[W_1^3]}{(\mathbb{E}[|W_1|^2])^{3/2}}\frac{1}{6\sqrt{2\pi n}}\right|\\ 
&\quad\le \frac{C}{6}\left[\frac{\mathbb{E}[|W_1|^4]}{\sigma^4n} + n^6\left(\sup_{|b| \ge \sigma^2/(12\beta_3)}|v(b)| + \frac{1}{2n}\right)^n\right].
\end{split}
\end{equation}
\end{proposition}
It is easy to see that $\kappa(t) \to 0$ as $t\to\infty$ and $\kappa(t) \to 2/3$ as $t\to-\infty$. Inequality~\eqref{eq:Delta-ReLU-each-t} in Proposition~\ref{prop:Delta-ReLU-bound} depends on the characteristic function $v(\cdot)$ of the random variable $W_1$. This can potentially be a sub-optimality of the bound. \cite{Borisov1996} proved a bound for $|\mathbb{E}[f((n\sigma^2)^{-1/2}\sum_{i=1}^n W_i)] - \mathbb{E}[f(Z)]|$ for twice differentiable function $f$ without any restriction on the characteristic function. We believe one can apply a similar technique to get a better bound for $\Delta_{n,\mathrm{ReLU}}(t)$, but we leave this for future work. 

Inequality~\eqref{eq:Delta-ReLU-integral} is useful in the context of Theorem~\ref{thm:integral-representation-ridge}, and in particular,~\eqref{eq:random-variable-ReLU-representation-assumption}. Furthermore, inequality~\eqref{eq:Delta-ReLU-integral} is also related to the bounds for ideal metrics. For example, 1-Wasserstein distance (or ideal metric of order 1) between two random variables $U$ and $V$ is given by
\[
d_1^{\mathrm{Wass}}(U, V) = \sup_{f\in\mathcal{F}_1}|\mathbb{E}[f(U)] - \mathbb{E}[f(V)]| = \int_{\mathbb{R}}|\mathbb{P}(U \le s) - \mathbb{P}(V \le s)|ds,
\]
where $\mathcal{F}_1 := \{f:\mathbb{R}\to\mathbb{R}\,|\,\sup_{x\neq y}|f(x) - f(y)|/|x - y| \le 1\}$ is the class of all $1$-Lipschitz functions.
Clearly, $\Delta_{n,\mathrm{ReLU}}(t) \le d_1^{\mathrm{Wass}}((n\sigma^2)^{-1/2}\sum_{i=1}^n W_i,\, Z)$ for all $t\in\mathbb{R}$. Bounds on $d_1^{\mathrm{Wass}}(\cdot, \cdot)$ for scaled averages of independent/dependent random variables are referred to as $L_1$ Berry--Esseen bounds~\citep{erickson1973l_p}. Although $L_1$ Berry--Esseen bounds are not as precise as~\eqref{eq:Delta-ReLU-integral}, they may suffice for the purpose of obtaining bounds for $\Delta_f(W_n, W)$ using Theorem 3.2 or~\eqref{eq:integral-representation-Radon-transform}. Optimal order $L_1$ Berry--Esseen bounds for the average of independent/dependent random variables can be found in~\cite{erickson1974l_1},~\cite{chen1986rate},~\cite{chen2004normal},~\cite{Goldstein2010},~\cite{van20141},~\cite{sunklodas2007normal}, and~\cite{Fan2020}. In particular, Theorems 3.2--3.4 of~\cite{bentkus2003new} imply bounds for ideal metrics of order $k \ge 1$ for sums of independent random variables, Theorems 1-2 of~\cite{sunklodas2007normal} imply bounds for the ideal metric of order 1 for sums of strongly mixing (or $\alpha$-mixing) random variables, Theorems 4-5 of~\cite{van20141} imply bounds for the ideal metric of order 1 for martingales. (Note that $L_1$ Berry--Esseen bounds can be obtained from non-uniform Berry--Esseen bounds.)
Inequality~\eqref{eq:Delta-ReLU-integral} suffices for using Theorem~\ref{thm:integral-representation-ridge}. Inequality~\eqref{eq:Delta-ReLU-integral} is related to ideal metrics of order 2; see Section 2.10 of~\cite{senatov2011normal}. The ideal metric of order 2 is given by
\begin{equation}\label{eq:ideal-metric-2}
\zeta_2(U, V) := \sup_{f\in\mathcal{F}_2}|\mathbb{E}[f(U)] - \mathbb{E}[f(V)]| = \int_{\mathbb{R}}\left|\mathbb{E}[(U - t)_+] - \mathbb{E}[(V - t)_+]\right|dt,
\end{equation}
where $\mathcal{F}_2 = \{f:\mathbb{R}\to\mathbb{R}\,|\,\|f^{(2)}\|_{\infty} \le 1\}$ is the class of all functions whose second derivative is uniformly bounded by 1.

What is interesting to note from Proposition~\ref{prop:Delta-ReLU-bound}~\eqref{eq:Delta-ReLU-each-t} is that if either $\mathbb{E}[W_1^3] = 0$ or $t = 0$, we get
\[
|\Delta_{n,\mathrm{ReLU}}(t)| ~\le~  C\left[\frac{\mathbb{E}[|W_1|^4]}{\sigma^4n} + n^6\left(\sup_{|b| \ge \sigma^2/(12\beta_3)}|v(b)| + \frac{1}{2n}\right)^n\right]\kappa(t),
\]
where the right hand side converges to zero at an $n^{-1}$ rate as $n\to\infty$. A simple implication is for the rate of convergence of moments. For example, note that $|x| = (x)_+ - (-x)_+$ and hence,
\[
\left|\mathbb{E}\left|\frac{1}{\sqrt{n\sigma^2}}\sum_{i=1}^n W_i\right| - \mathbb{E}[|Z|]\right| ~\le~ 2C\left[\frac{\mathbb{E}[|W_1|^4]}{\sigma^4n} + n^6\left(\sup_{|b| \ge \sigma^2/(12\beta_3)}|v(b)| + \frac{1}{2n}\right)^n\right]\kappa(0).
\]
Lemma 6.2 of~\cite{kainen2010integral} and Example 4.8(2) of~\cite{wojtowytsch2020representation} prove that 
\[
\|x\|_2 ~=~ c_d\int_{S^{d-1}} (\langle a, x\rangle)_+\pi^0(a)da,\quad\mbox{for all}\quad x\in\mathbb{R}^d,
\]
where $\pi^0(\cdot)$ is the uniform measure on the unit sphere $(S^{d-1}=\{\theta\in\mathbb{R}^d:\,\|\theta\|_2 = 1\})$ and $$c_d = \left(\int_{S^{d-1}} (\langle e_1, w\rangle)_+\pi^0(w)dw\right)^{-1}~\asymp~ 2\sqrt{\pi d}.$$ 
Therefore, for mean zero independent and identically distributed random vectors $X_1, \ldots, X_n\in\mathbb{R}^d$, we have
\begin{align*}
&\left|\mathbb{E}\left[\left\|\frac{1}{\sqrt{n}}\sum_{i=1}^n X_i\right\|_2\right] - \mathbb{E}[\|Z'\|_2]\right|\\ 
&\quad\le 2C\kappa(0)c_d\int_{S^{d-1}} \left[\frac{\mathbb{E}[|\langle a, X_1\rangle|^4]}{(\mathbb{E}[|\langle a, W_1\rangle|^2])^{3/2}n} + n^6(\mathbb{E}[|\langle a, X_1\rangle|^2])^{1/2}\left(\sup_{|b| \ge \gamma_a}|v_a(b)| + \frac{1}{2n}\right)^n\right]\pi^0(a)da, 
\end{align*}
where $v_a(b) = \mathbb{E}[e^{i\langle a, X_1\rangle}]$ and $\gamma_a = \mathbb{E}[|\langle a, X_1\rangle|^2]/(\mathbb{E}[|\langle a, X_1\rangle|^3])^{3/2}$.
One can simply bound the integral on the right hand side by the supremum over all $a\in S^{d-1}$. If the random vector $X_1\in\mathbb{R}^d$ satisfies $L_4$--$L_2$ moment equivalence (i.e., there exists a constant $L$ such that $\mathbb{E}[|\langle a, X_1\rangle|^4] \le L(\mathbb{E}[|\langle a, X_1\rangle|^2])^2$ for all $a\in S^{d-1}$), then $1/\gamma_a \le L^3\|\Sigma\|_{op}^{1/2}$ which implies that $$\sup_{|b| \ge \gamma_a}|v_a(b)| \le \sup_{|b| \ge L^{-3}\|\Sigma\|_{op}^{-1/2}}|v_a(b)|,$$ and hence, we get
\begin{equation}\label{eq:Euclidean-norm-expectations}
\begin{split}
&\left|\mathbb{E}\left[\left\|\frac{1}{\sqrt{n}}\sum_{i=1}^n X_i\right\|_2\right] - \mathbb{E}[\|Z'\|_2]\right|\\ 
&\quad\le 2C\kappa(0)c_d\|\Sigma\|_{op}^{1/2}\left[\frac{L}{n} + n^6\int_{S^{d-1}}\left(\sup_{|b| \ge L^{-3}\|\Sigma\|_{op}^{-1/2}}|v_a(b)| + \frac{1}{2n}\right)^n\pi^0(a)da\right].
\end{split}
\end{equation}
The assumption of $L_4$--$L_2$ moment equivalence is standard in robust covariance estimation as well as small ball property; see~\cite{oliveira2016lower} and~\cite{mendelson2020robust}. One important class of distributions that satisfy $L_4$--$L_2$ moment equivalence assumption is the class of log-concave distributions; see Remark 2.20 of~\cite{patil2022mitigating}.

If the characteristic function $v_a(b)$ is bounded away from 1 for almost all $a\in S^{d-1}$, then the second term on the right hand side of~\eqref{eq:Euclidean-norm-expectations} converges to zero exponentially in $n$ and hence, we get that the difference of the expectations of $\|\cdot\|_2$-norms is of order $c_d/n\asymp d^{1/2}/n$ as $n\to\infty$. There are two interesting features of this result: (1) the rate of convergence is $n^{-1},$ and (2) the right hand side converges to zero even with increasing dimension as long as $d = o(n^2)$. The rate of convergence of $n^{-1}$ for the difference of probabilities of Euclidean balls was proved in~\cite{Bentkus1996} and~\cite{gotze2014explicit} for $d \ge 5$. The dependence on dimension $d$, however, is not known. To the best of our knowledge, inequality similar to~\eqref{eq:Euclidean-norm-expectations} is unknown in the literature. It may be worth noting here that $f(x) = \|x\|_2$ is not a smooth function; it is not differentiable at $x = 0$.
\section{Discussion}\label{sec:discussion}
We have discussed ways to bound the differences of expectations of averages of random variables using Berry--Esseen results for Euclidean balls or for univariate random variables. We have also summarized several integral representations of functions that can yield bounds that have none to minimal dependence on the dimension of the underlying random variables. Two interesting aspects of our bounds are that (1) they are equally applicable to both independent and dependent random variables and (2) they can be used with any limiting distribution (or any infinitely divisible distribution). All the results presented in this paper require some sort of smoothness on the functions expressed either by weighted integrability of the Fourier transform, or explicitly by the differentiability requirement. Obtaining bounds for $\Delta_f(W_n, W)$ for arbitrary Borel measurable function $f$ is usually done by smoothing with a mollifier and it seems plausible that one can apply the results of this paper in that context.
\paragraph{Acknowledgments.} This work is partially supported by NSF DMS–2113611.
\bibliography{references}

\begin{thebibliography}{}

\bibitem[Abdeljawad and Grohs, 2022]{abdeljawad2022integral}
Abdeljawad, A. and Grohs, P. (2022).
\newblock Integral representations of shallow neural network with rectified
  power unit activation function.
\newblock {\em Neural Networks}, 155:536--550.

\bibitem[Banerjee and Kuchibhotla, 2023]{banerjee2023central}
Banerjee, A. and Kuchibhotla, A.~K. (2023).
\newblock Central limit theorems and approximation theory: Part {I}.
\newblock {\em arXiv preprint arXiv:2306.05947}.

\bibitem[Barron, 1993]{barron1993universal}
Barron, A.~R. (1993).
\newblock Universal approximation bounds for superpositions of a sigmoidal
  function.
\newblock {\em IEEE Transactions on Information theory}, 39(3):930--945.

\bibitem[Bentkus, 2003a]{bentkus2003new}
Bentkus, V. (2003a).
\newblock A new method for approximations in probability and operator theories.
\newblock {\em Lithuanian Mathematical Journal}, 43:367--388.

\bibitem[Bentkus, 2003b]{bentkus2003dependence}
Bentkus, V. (2003b).
\newblock On the dependence of the berry--esseen bound on dimension.
\newblock {\em Journal of Statistical Planning and Inference}, 113(2):385--402.

\bibitem[Bentkus, 2004]{bentkus2005lyapunov}
Bentkus, V. (2004).
\newblock A {L}yapunov type bound in {${\bf R}^d$}.
\newblock {\em Teor. Veroyatn. Primen.}, 49(2):400--410.

\bibitem[Bentkus and G{\"o}tze, 1993]{bentkus1993smoothness}
Bentkus, V. and G{\"o}tze, F. (1993).
\newblock On smoothness conditions and convergence rates in the clt in banach
  spaces.
\newblock {\em Probability theory and related fields}, 96(2):137--151.

\bibitem[Bentkus and G\"{o}tze, 1996]{Bentkus1996}
Bentkus, V. and G\"{o}tze, F. (1996).
\newblock Optimal rates of convergence in the {CLT} for quadratic forms.
\newblock {\em Ann. Probab.}, 24(1):466--490.

\bibitem[Bhattacharya and Rao, 2010]{bhattacharya2010normal}
Bhattacharya, R.~N. and Rao, R.~R. (2010).
\newblock {\em Normal approximation and asymptotic expansions}.
\newblock SIAM.

\bibitem[Bogatyrev et~al., 2006]{bogatyrev2006non}
Bogatyrev, S., G{\"o}tze, F., and Ulyanov, V. (2006).
\newblock Non-uniform bounds for short asymptotic expansions in the {CLT} for
  balls in a {H}ilbert space.
\newblock {\em Journal of multivariate analysis}, 97(9):2041--2056.

\bibitem[Bong et~al., 2022]{bong2022high}
Bong, H., Kuchibhotla, A.~K., and Rinaldo, A. (2022).
\newblock High-dimensional berry-esseen bound for $ m $-dependent random
  samples.
\newblock {\em arXiv preprint arXiv:2212.05355}.

\bibitem[Borisov and Skilyagina, 1996]{Borisov1996}
Borisov, I.~S. and Skilyagina, G.~I. (1996).
\newblock On the asymptotic expansion of the moments of smooth functions in the
  central limit theorem.
\newblock {\em Sibirsk. Mat. Zh.}, 37(3):519--525, i.

\bibitem[Chen, 1986]{chen1986rate}
Chen, L.~H. (1986).
\newblock The rate of convergence in a central limit theorem for dependent
  random variables with arbitrary index set.

\bibitem[Chen and Shao, 2004]{chen2004normal}
Chen, L.~H. and Shao, Q.-M. (2004).
\newblock Normal approximation under local dependence.
\newblock {\em Annals of probability}, 32(3A):1985--2028.

\bibitem[Chernozhukov et~al., 2023]{chernozhukov2023nearly}
Chernozhukov, V., Chetverikov, D., and Koike, Y. (2023).
\newblock Nearly optimal central limit theorem and bootstrap approximations in
  high dimensions.
\newblock {\em The Annals of Applied Probability}, 33(3):2374--2425.

\bibitem[De~Vito et~al., 2021]{de2021reproducing}
De~Vito, E., M{\"u}cke, N., and Rosasco, L. (2021).
\newblock Reproducing kernel hilbert spaces on manifolds: Sobolev and diffusion
  spaces.
\newblock {\em Analysis and Applications}, 19(03):363--396.

\bibitem[Diewert et~al., 1981]{diewert1981nine}
Diewert, W.~E., Avriel, M., and Zang, I. (1981).
\newblock Nine kinds of quasiconcavity and concavity.
\newblock {\em Journal of Economic Theory}, 25(3):397--420.

\bibitem[Erickson, 1973]{erickson1973l_p}
Erickson, R. (1973).
\newblock On an l\_p version of the berry-esseen theorem for independent and
  m-dependent variables.
\newblock {\em The Annals of Probability}, pages 497--503.

\bibitem[Erickson, 1974]{erickson1974l_1}
Erickson, R. (1974).
\newblock L\_1 bounds for asymptotic normality of m-dependent sums using
  stein's technique.
\newblock {\em The Annals of Probability}, pages 522--529.

\bibitem[Fan and Ma, 2020]{Fan2020}
Fan, X. and Ma, X. (2020).
\newblock On the {W}asserstein distance for a martingale central limit theorem.
\newblock {\em Statist. Probab. Lett.}, 167:108892, 6.

\bibitem[Fang et~al., 2023]{fang2023high}
Fang, X., Koike, Y., Liu, S.-H., and Zhao, Y.-K. (2023).
\newblock High-dimensional central limit theorems by stein's method in the
  degenerate case.
\newblock {\em arXiv preprint arXiv:2305.17365}.

\bibitem[Fomin, 1983]{fomin1983estimate}
Fomin, S. (1983).
\newblock An estimate of the rate of convergence in the multi-dimensional
  central limit theorem.
\newblock {\em Theory of Probability \& Its Applications}, 27(2):365--368.

\bibitem[Funahashi, 1989]{funahashi1989approximate}
Funahashi, K.-I. (1989).
\newblock On the approximate realization of continuous mappings by neural
  networks.
\newblock {\em Neural networks}, 2(3):183--192.

\bibitem[Gin{\'e} and Nickl, 2021]{gine2021mathematical}
Gin{\'e}, E. and Nickl, R. (2021).
\newblock {\em Mathematical foundations of infinite-dimensional statistical
  models}.
\newblock Cambridge university press.

\bibitem[Girosi, 1994]{girosi1994regularization}
Girosi, F. (1994).
\newblock Regularization theory, radial basis functions, and networks.
\newblock In {\em From statistics to neural networks: Theory and pattern
  recognition applications}, pages 166--187. Springer.

\bibitem[Girosi and Anzellotti, 1992]{girosi1992convergence}
Girosi, F. and Anzellotti, G. (1992).
\newblock Convergence rates of approximation by translates.
\newblock Technical report, M.I.T. AI Memo. No.1288, MIT, Cambridge, MA.

\bibitem[Goldstein, 2010]{Goldstein2010}
Goldstein, L. (2010).
\newblock Bounds on the constant in the mean central limit theorem.
\newblock {\em Ann. Probab.}, 38(4):1672--1689.

\bibitem[Gotze, 1986]{gotze1986rate}
Gotze, F. (1986).
\newblock On the rate of convergence in the central limit theorem in banach
  spaces.
\newblock {\em The Annals of Probability}, pages 922--942.

\bibitem[G{\"o}tze and Zaitsev, 2014]{gotze2014explicit}
G{\"o}tze, F. and Zaitsev, A.~Y. (2014).
\newblock Explicit rates of approximation in the clt for quadratic forms.
\newblock {\em The Annals of Probability}, 42(1):354--397.

\bibitem[Heinrich, 1985]{heinrich1985non}
Heinrich, L. (1985).
\newblock Non-uniform estimates, moderate and large deviations in the central
  limit theorem for m-dependent random variables.
\newblock {\em Mathematische Nachrichten}, 121(1):107--121.

\bibitem[Hipp, 1979]{hipp1979convergence}
Hipp, C. (1979).
\newblock Convergence rates in the central limit theorem for stationary mixing
  sequences of random vectors.
\newblock {\em Journal of Multivariate Analysis}, 9(4):560--578.

\bibitem[Irie and Miyake, 1988]{irie1988capabilities}
Irie, B. and Miyake, S. (1988).
\newblock Capabilities of three-layered perceptrons.
\newblock In {\em ICNN}, pages 641--648.

\bibitem[Jirak, 2015]{jirak2015non}
Jirak, M. (2015).
\newblock On non-uniform berry--esseen bounds for time series.
\newblock {\em Probab. Math. Statist}, 35(1):1--14.

\bibitem[Kainen et~al., 2010]{kainen2010integral}
Kainen, P.~C., K{\u u}rkov{\'a}, V., and Vogt, A. (2010).
\newblock Integral combinations of heavisides.
\newblock {\em Mathematische Nachrichten}, 283(6):854--878.

\bibitem[Klusowski and Barron, 2018]{klusowski2018approximation}
Klusowski, J.~M. and Barron, A.~R. (2018).
\newblock Approximation by combinations of relu and squared relu ridge
  functions with $\ell^{} 1$ and $\ell^{} 0$ controls.
\newblock {\em IEEE Transactions on Information Theory}, 64(12):7649--7656.

\bibitem[Makovoz, 1998]{makovoz1998uniform}
Makovoz, Y. (1998).
\newblock Uniform approximation by neural networks.
\newblock {\em Journal of Approximation Theory}, 95(2):215--228.

\bibitem[Mendelson and Zhivotovskiy, 2020]{mendelson2020robust}
Mendelson, S. and Zhivotovskiy, N. (2020).
\newblock Robust covariance estimation under $l_4$--$l_2$ norm equivalence.
\newblock {\em The Annals of Statistics}, 48(3):1648--1664.

\bibitem[Michel, 1978]{michel1978remark}
Michel, R. (1978).
\newblock A remark on nonuniform estimates in the central limit theorem for
  sums of independent random vectors.
\newblock {\em Sankhy{\=a}: The Indian Journal of Statistics, Series A}, pages
  388--392.

\bibitem[Minh, 2010]{minh2010some}
Minh, H.~Q. (2010).
\newblock Some properties of gaussian reproducing kernel hilbert spaces and
  their implications for function approximation and learning theory.
\newblock {\em Constructive Approximation}, 32(2):307--338.

\bibitem[Novak et~al., 2018]{novak2018reproducing}
Novak, E., Ullrich, M., Wo{\'z}niakowski, H., and Zhang, S. (2018).
\newblock Reproducing kernels of sobolev spaces on $\mathbb{R}^d$ and
  applications to embedding constants and tractability.
\newblock {\em Analysis and Applications}, 16(05):693--715.

\bibitem[Oliveira, 2016]{oliveira2016lower}
Oliveira, R.~I. (2016).
\newblock The lower tail of random quadratic forms with applications to
  ordinary least squares.
\newblock {\em Probability Theory and Related Fields}, 166:1175--1194.

\bibitem[Patil et~al., 2022]{patil2022mitigating}
Patil, P., Kuchibhotla, A.~K., Wei, Y., and Rinaldo, A. (2022).
\newblock Mitigating multiple descents: A model-agnostic framework for risk
  monotonization.
\newblock {\em arXiv preprint arXiv:2205.12937}.

\bibitem[Paulauskas, 1976]{paulauskas1976convergence}
Paulauskas, V. (1976).
\newblock Convergence of some functionals of sums of independent random
  variables in a banach space.
\newblock {\em Lithuanian Mathematical Journal}, 16(3):385--399.

\bibitem[Paulauskas and Ra{\v c}kauskas, 1989]{paulauskas_Rackauskas_2012}
Paulauskas, V. and Ra{\v c}kauskas, A. (1989).
\newblock {\em Approximation theory in the central limit theorem}, volume~32 of
  {\em Mathematics and its Applications (Soviet Series)}.
\newblock Kluwer Academic Publishers Group, Dordrecht.
\newblock Exact results in Banach spaces, Translated from the Russian by B.
  Svecevi\v cius and Paulauskas.

\bibitem[Petrosyan et~al., 2020]{petrosyan2020neural}
Petrosyan, A., Dereventsov, A., and Webster, C.~G. (2020).
\newblock Neural network integral representations with the relu activation
  function.
\newblock In {\em Mathematical and Scientific Machine Learning}, pages
  128--143. PMLR.

\bibitem[Petrov, 1995]{Petrov1995}
Petrov, V.~V. (1995).
\newblock {\em Limit theorems of probability theory}, volume~4 of {\em Oxford
  Studies in Probability}.
\newblock The Clarendon Press, Oxford University Press, New York.
\newblock Sequences of independent random variables, Oxford Science
  Publications.

\bibitem[Rachev and Yukich, 1989]{rachev1989rates}
Rachev, S. and Yukich, J. (1989).
\newblock Rates for the clt via new ideal metrics.
\newblock {\em The Annals of Probability}, pages 775--788.

\bibitem[Rai{\v{c}}, 2019]{raivc2019multivariate}
Rai{\v{c}}, M. (2019).
\newblock A multivariate berry--esseen theorem with explicit constants.
\newblock {\em Bernoulli}, 25(4A):2824--2853.

\bibitem[Rotar’, 1970]{rotar1970non}
Rotar’, V.~I. (1970).
\newblock A non-uniform estimate for the convergence speed in the
  multi-dimensional central theorem.
\newblock {\em Theory of Probability \& Its Applications}, 15(4):630--648.

\bibitem[Sazonov, 1974]{sazonov1974new}
Sazonov, V. (1974).
\newblock A new general estimate of the rate of convergence in the central
  limit theorem in rk.
\newblock {\em Proceedings of the National Academy of Sciences},
  71(1):118--121.

\bibitem[Schaback, 2007]{schaback2007kernel}
Schaback, R. (2007).
\newblock Kernel-based meshless methods.
\newblock {\em Lecture Notes for Taught Course in Approximation Theory.
  Georg-August-Universit{\"a}t G{\"o}ttingen}.

\bibitem[Senatov, 1982]{senatov1982some}
Senatov, V. (1982).
\newblock Some non-uniform estimates of the speed of convergence in the
  multi-dimensional central limit theorem.
\newblock {\em Theory of Probability \& Its Applications}, 26(4):657--669.

\bibitem[Senatov, 1992]{senatov1992some}
Senatov, V.~V. (1992).
\newblock Some remarks on estimating the rate of convergence in the central
  limit theorem in the hilbert space.
\newblock {\em Theory of Probability \& Its Applications}, 36(2):401--405.

\bibitem[Senatov, 1993]{senatov1993estimates}
Senatov, V.~V. (1993).
\newblock On estimates of the rate of convergence in the central limit theorem
  in multidimensional spaces.
\newblock {\em Theory of Probability \& Its Applications}, 37(4):703--706.

\bibitem[Senatov, 2011]{senatov2011normal}
Senatov, V.~V. (2011).
\newblock {\em Normal approximation: new results, methods and problems}.
\newblock Walter de Gruyter.

\bibitem[Siegel and Xu, 2020]{siegel2020approximation}
Siegel, J.~W. and Xu, J. (2020).
\newblock Approximation rates for neural networks with general activation
  functions.
\newblock {\em Neural Networks}, 128:313--321.

\bibitem[Steinwart and Christmann, 2008]{steinwart2008support}
Steinwart, I. and Christmann, A. (2008).
\newblock {\em Support vector machines}.
\newblock Springer Science \& Business Media.

\bibitem[Stuetzle and Mittal, 2006]{stuetzle2006some}
Stuetzle, W. and Mittal, Y. (2006).
\newblock Some comments on the asymptotic behavior of robust smoothers.
\newblock In {\em Smoothing Techniques for Curve Estimation: Proceedings of a
  Workshop Held in Heidelberg, April 2--4, 1979}, pages 191--195. Springer.

\bibitem[Sunklodas, 2007]{sunklodas2007normal}
Sunklodas, J. (2007).
\newblock On normal approximation for strongly mixing random variables.
\newblock {\em Acta Applicandae Mathematicae}, 97(1-3):251--260.

\bibitem[Tikhomirov, 1994]{tikhomirov1994accuracy}
Tikhomirov, A.~N. (1994).
\newblock On the accuracy of normal approximation of the probability of hitting
  a ball of sums of weakly dependent hilbert space valued random variables ii.
\newblock {\em Theory of Probability \& Its Applications}, 38(1):80--94.

\bibitem[Ulyanov, 1986]{ulyanov1986normal}
Ulyanov, V. (1986).
\newblock Normal approximation for sums of nonidentically distributed random
  variables in hilbert spaces.
\newblock {\em Acta Scientiarum Mathematicarum}, 50(3-4):411--419.

\bibitem[Van~Dung et~al., 2014]{van20141}
Van~Dung, L., Son, T.~C., and Tien, N.~D. (2014).
\newblock L 1 bounds for some martingale central limit theorems.
\newblock {\em Lithuanian Mathematical Journal}, 54(1):48--60.

\bibitem[Wainwright, 2019]{wainwright2019high}
Wainwright, M.~J. (2019).
\newblock {\em High-dimensional statistics: A non-asymptotic viewpoint},
  volume~48.
\newblock Cambridge university press.

\bibitem[Weinan and Wojtowytsch, 2022]{wojtowytsch2020representation}
Weinan, E. and Wojtowytsch, S. (2022).
\newblock Representation formulas and pointwise properties for barron
  functions.
\newblock {\em Calculus of Variations and Partial Differential Equations},
  61(2):46.

\bibitem[Yaroslavtseva, 2008]{yaroslavtseva2008non}
Yaroslavtseva, L. (2008).
\newblock {\em Non-classical error bounds in the central limit theorem}.
\newblock PhD thesis.
\newblock (Doctoral dissertation, Otto-von-Guericke-Universität Magdeburg,
  Universit{\"a}tsbibliothek).

\bibitem[Yurinskii, 1983]{yurinskii1983accuracy}
Yurinskii, V. (1983).
\newblock On the accuracy of normal approximation of the probability of hitting
  a ball.
\newblock {\em Theory of Probability \& Its Applications}, 27(2):280--289.

\bibitem[Zlobec, 2006]{zlobec2006characterization}
Zlobec, S. (2006).
\newblock Characterization of convexifiable functions.
\newblock {\em Optimization}, 55(3):251--261.

\end{thebibliography}
\bibliographystyle{apalike}


\newpage
\appendix
\section{Proof of Proposition~\ref{prop:signed-measure-bound}}\label{sec:proof-signed-measure-bound}
By Tonelli's theorem, we have
\[
\mathbb{E}[g(W_n)] = \int_{\Lambda} \mathbb{E}[f_{\lambda}(W_n)]\Lambda(d\lambda)\quad\mbox{and}\quad \mathbb{E}[g(W)] = \int \mathbb{E}[f_{\lambda}(W)]\mu(d\lambda).
\]
This implies that 
\begin{align*}
\Delta_g(W_n, W) &= \left|\int_{\Lambda} \left\{\mathbb{E}[f_{\lambda}(W_n)] - \mathbb{E}[f_{\lambda}(W)]\right\}\mu(d\lambda)\right|\\
&\le \int_{\Lambda} \left|\mathbb{E}[f_{\lambda}(W_n)] - \mathbb{E}[f_{\lambda}(W)]\right||\mu|(d\lambda)\\
&= \int_{\Lambda} \Delta_{f_{\lambda}}(W_n, W)|\mu|(d\lambda).
\end{align*}
This completes the proof.
\section{Proof of Norm Ball Result}\label{sec:proof-of-norm-ball-result}
Set
\[
C_1 = \frac{C(2\pi)^{-d/2}\beta_3}{n^{1/2}h^d},\quad a = \frac{1}{\sigma_1\cdots\sigma_6},\quad\mbox{and}\quad b = \frac{1}{\sigma^3} + \frac{1}{(\sigma_1\cdots\sigma_6)^{1/2}}
\]
From the bound on $\Delta_{f_h}(W_n, W)$, the bound we need to compute is of the form
\[
C_1\int_{\mathbb{R}^d\times[0, 1]} \left(a(2h\log(1/t))^{3/2} + b\right)\left[\mathbf{1}\{\|y\|_2 < \sqrt{8h\log(1/t)}\} + \frac{\mathbf{1}\{\|y\|_2 \ge \sqrt{8h\log(1/t)}\}}{1 + \|y\|_2^3/(8\sigma^3)}\right]f(y)dtdy.
\]
We now compute the integral first with respect to $t$. By H{\"o}lder's inequality (if $1/p + 1/q = 1$),
\begin{equation*}
\begin{split}
\int_0^1 &(2h\log(1/t))^{3/2}\mathbf{1}\{\|y\|_2 < \sqrt{8h\log(1/t)}\}dt\\ 
&\le \left(\int_0^1 (2h\log(1/t))^{3p/2}dt\right)^{1/p}\left(\int_0^1 \mathbf{1}\{\|y\|_2 < \sqrt{8h\log(1/t)}\}dt\right)^{1/q} \\
&\le (2h)^{3/2}\left(\Gamma\left(\frac{3p}{2} + 1\right)\right)^{1/p}\left(\int_0^{\exp(-\|y\|_2^2/(8h))} dt\right)^{1/q}\\
&\le (3/2)^{5/2}e^{1/e}p^{3/2}(2h)^{3/2}\exp\left(-\frac{\|y\|_2^2}{8hq}\right),
\end{split}
\end{equation*}
where the last inequality follows from the relations $\Gamma(x + 1) = x\Gamma(x)$ and $\Gamma(x) \le x^x$ for all $x > 0$.
One could potentially optimize over $p\ge1$ and $q = (1 - 1/p)^{-1}$. But we simply take $p = q = 2$, for brevity and obtain
\begin{equation}\label{eq:product-3/2-indicator}
\int_0^1 (2h\log(1/t))^{3/2}\mathbf{1}\{\|y\|_2 < \sqrt{8h\log(1/t)}\}dt ~\le~ Ch^{3/2}\exp\left(-\frac{\|y\|_2^2}{16h}\right),
\end{equation}
for some absolute constant $C = (3/2)^{5/2}e^{1/e}2^{3}$.
Similarly, 
\begin{equation}\label{eq:constant-indicator-less-and-greater}
\begin{split}
\int_{0}^1 \mathbf{1}\{\|y\|_2 < \sqrt{8h\log(1/t)}\}dt ~&=~ \exp\left(-\frac{\|y\|_2^2}{8h}\right),\\
\int_0^1 \mathbf{1}\{\|y\|_2 \ge \sqrt{8h\log(1/t)}\}dt ~&\le~ 1,
\end{split}
\end{equation}
and
\begin{equation}\label{eq:product-3/2-greater}
\begin{split}
\int_0^1 (2h\log(1/t))^{3/2}\mathbf{1}\{\|y\|_2 \ge \sqrt{8h\log(1/t)}\}dt \le \int_0^1 (2h\log(1/t))^{3/2} = (2h)^{3/2}\Gamma(5/2).
\end{split}
\end{equation}
Combining~\eqref{eq:product-3/2-indicator},~\eqref{eq:constant-indicator-less-and-greater}, and~\eqref{eq:product-3/2-greater}, we conclude that (for some absolute constant $C$)
\begin{align*}
&\Delta_{f_h}(W_n, W)\\ 
&\le CC_1\int_{\mathbb{R}^d} \left[ah^{3/2}\exp\left(-\frac{\|y\|_2^2}{16h}\right) + b\exp\left(-\frac{\|y\|_2^2}{8h}\right) + \frac{(ah^{3/2} + b)}{1 + \|y\|_2^3/(8\sigma^3)}\right]f(y)dy\\
&\le CC_1(ah^{3/2} + b)\int_{\mathbb{R}^d}\left[\exp\left(-\frac{\|y\|_2^2}{16h}\right) + \frac{1}{1 + \|y\|_2^3/(8\sigma^3)}\right]f(y)dy\\
&\le \frac{C\beta_3}{n^{1/2}(2\pi)^{d/2}}\left(\frac{h^{3/2}}{\sigma_1\cdots\sigma_6} + \frac{1}{\sigma^3} + \frac{1}{(\sigma_1\cdots\sigma_6)^{1/2}}\right)\left(\int_{\mathbb{R}^d} \left[\frac{1}{1 + h^3\|s\|_2^3/(8\sigma^3)} +  \exp\left(-\frac{\|s\|_2^2}{16}\right)\right]f(hs)ds\right)
\end{align*}
\section{Proof of Theorem~\ref{thm:integral-representation-ridge}}\label{sec:proof-of-integral-representation}
The main idea behind the proof is to note that 
\[
f(x) - f(0) - f'(0)[x] = \int_{\mathbb{R}^d} (e^{-i\langle \omega, x\rangle} - i\langle \omega, x\rangle - 1)\widecheck{f}(\omega)d\omega,
\]
where $f'(0)[x] = (d/dt)f(xt)\big|_{t = 0}$ is the directional derivative of $f$ at $0$ in the direction of $x$.
We first represent $z\mapsto e^{iz} - iz - 1$ in terms of the ReLU function so that $f(x)$ can be decomposed into a linear function of $x$ and an integral of ReLU functions.
Note the identity
\[
-\int_0^{\infty} \left[(z - u)_+e^{iu} + (-z-u)e^{-iu}\right]du = e^{iz} - iz - 1,\quad\mbox{for all}\quad z\in\mathbb{R}.
\]
In particular, taking $z = \langle\omega, x\rangle$, we get
\begin{equation}\label{eq:penultimate-display}
-\int_0^{\infty} \left[(\langle\omega, x\rangle - u)_+e^{iu} + (-\langle\omega, x\rangle - u)_+e^{-iu}\right]du = e^{i\langle \omega, x\rangle} - i\langle \omega, x\rangle - 1.    
\end{equation}
Because $|e^{iz} - iz - 1| \le \min\{2|z|,\,|z|^2/2\}$, we get that $\int_{\mathbb{R}^d}(e^{i\langle \omega, x\rangle} - i\langle \omega, x\rangle - 1)\widecheck{f}(\omega)d\omega$ exists as long as $\int_{\mathbb{R}^d} \min\{2|\langle\omega, x\rangle|, |\langle \omega, x\rangle|^2/2\}|\widecheck{f}(\omega)|d\omega < \infty$, which holds by~\eqref{eq:ReLU-representation-assumption}. Integrating both sides of~\eqref{eq:penultimate-display} with respect to $\widecheck{f}(\omega)$ over $\omega\in\mathbb{R}^d$, we get
\[
-\int_{\mathbb{R}^d}\int_0^{\infty} \left[(\langle\omega,x\rangle - u)_+e^{iu} + (-\langle\omega, x\rangle-u)_+e^{-iu}\right]\widecheck{f}(\omega)dud\omega = f(x) - f(0) - f'(0)[x].
\]
Because the right hand side is real-valued, the real part of the left hand side is the same as the right hand side as well. Therefore, writing $\widecheck{f}(\omega) = |\widecheck{f}(\omega)|e^{ib(\omega)}$ for some $b(\omega)$, we conclude that
\begin{align*}
f(x) &- f(0) - f'(0)[x]\\
&=-\int_{\mathbb{R}^d}\int_{0}^{\infty} \left[(\langle \omega, x\rangle - u)_+\cos(u + b(\omega)) + (-\langle\omega, x\rangle - u)_+\cos(u - b(\omega))\right]|\widecheck{f}(\omega)|dud\omega.
\end{align*}
The bound on $\Delta_f(W_n, W)$ follows similarly under~\eqref{eq:random-variable-ReLU-representation-assumption}.
\section{Proof of Theorem~\ref{thm:integral-representation-general-activation}}\label{sec:proof-of-integral-representation-general-activation}
The proof is given in~\cite{irie1988capabilities} and~\cite{siegel2020approximation}. We repeat it here for completeness. By definition of Fourier transform for functions on $\mathbb{R}$, we have
\[
0\neq \widecheck{h}(a) = \frac{1}{2\pi}\int_{\mathbb{R}} h(t)e^{-iat}dt.
\]
Applying the change of variable $t\to \langle \omega, x\rangle + u$ (for some fixed $\omega$ and $x$), we get
\[
\widecheck{h}(a) = \frac{1}{2\pi}\int_{\mathbb{R}} h(\langle\omega,x\rangle + u)e^{-ia(\langle\omega,x\rangle + u)}du.
\]
Rearranging, we conclude
\[
e^{ia\langle\omega,x\rangle} = \frac{1}{2\pi|\widecheck{h}(a)|}\int_{\mathbb{R}} h(\langle\omega,x\rangle + u)e^{-i(au + c_h)}du.
\]
This implies
\[
e^{i\langle\omega,x\rangle} = \frac{1}{2\pi|\widecheck{h}(a)|}\int_{\mathbb{R}} h(a^{-1}\langle\omega,x\rangle + u)e^{-i(au + c_h)}du.
\]
Because $f(x) = \int_{\mathbb{R}^d} e^{i\langle\omega,x\rangle}\widecheck{f}(\omega)d\omega$, we obtain
\[
f(x) = \int_{\mathbb{R}^d}\int_{\mathbb{R}} \frac{1}{2\pi|\widecheck{h}(a)|}h(a^{-1}\langle\omega,x\rangle + u)e^{-i(au + c_h)}\widecheck{f}(\omega)dud\omega.
\]
Because the left hand side is real-valued, we can take the real part of the right hand side to obtain the result.
\section{Proof of Proposition~\ref{prop:Delta-ReLU-bound}}\label{sec:proof-of-prop-Delta-ReLU-bound}
The proof follows from the following equalities.
Note that
\[
\int_0^{\infty} (1 - (t + s)^2)\frac{e^{-(t + s)^2/2}}{\sqrt{2\pi}}ds = \int_t^{\infty} (1 - x^2)\frac{e^{-x^2/2}}{\sqrt{2\pi}}dx = x\frac{e^{-x^2/2}}{\sqrt{2\pi}}\big|_t^{\infty} = -t\frac{e^{-t^2/2}}{\sqrt{2\pi}}.
\]
Moreover, we have
\begin{align*}
\int_0^{\infty} \frac{1}{(1 + |t + s|)^4}ds ~=~ \begin{cases}
(3(1 + t)^3)^{-1}, &\mbox{if }t \ge 0,\\
2/3 - (3(1-t)^3)^{-1}, &\mbox{if }t < 0.
\end{cases}
\end{align*}
\end{document}